\documentclass[11pt,a4paper]{article}
\usepackage[novbox]{pdfsync}
\usepackage{amsfonts,amsmath,amssymb,mathrsfs,url,
esint,subcaption,
comment}
\usepackage{fancyhdr}
\usepackage{float,multirow,lipsum}
\usepackage{array,listings,verbatim}
\usepackage{tabularx}
\usepackage{nccmath}
\usepackage[utf8]{inputenc}
\usepackage{todonotes}

\usepackage{graphicx}
\graphicspath{{figures/}}
\graphicspath{{Notebook/}}
\usepackage[normalem]{ulem}

\usepackage{lineno}
\usepackage{color}
\usepackage{textcomp}
\usepackage{geometry}
\usepackage{enumitem}

\usepackage{graphicx}
\usepackage{epsfig}
\usepackage{wasysym}
\usepackage{empheq}
\usepackage{url}
\usepackage{tikz}
\usetikzlibrary{positioning}
\usetikzlibrary{arrows}
\usetikzlibrary{arrows.meta}
\usetikzlibrary{calc}
\newtheorem{theorem}{Theorem}
\newtheorem{remark}{Remark}
\newtheorem{proposition}{Proposition}
\newtheorem{definition}{Definition}
\newtheorem{corollary}{Corollary}
\newtheorem{example}{Example}
\newtheorem{lemma}{Lemma}
\newtheorem{assumption}{Assumption}
\newtheorem{claim}{Claim}

 \def\eeD{\end{definition}} \def\beD{\begin{definition}}
\def\beR{\begin{remark}} \def\eeR{\end{remark}}
\def\beL{\begin{lemma}} \def\eeL{\end{lemma}}
\def\beC{\begin{corollary}
  }\def\eeC{\end{corollary}}
  \def\beT{\begin{theorem}}\def\eeT{\end{theorem}}
  \def\beP{\begin{proposition}} \def\eeP{\end{proposition}}
\def\beXa{\begin{example}} \def\eeXa{\end{example}}
\def\beA{\begin{assumption}} \def\eeA{\end{assumption}}

\newtheorem{prop}{Proposition}

\newtheorem{ass}{Assumption}
\def\beAs{\begin{ass}
  }
\def\eeAs{\end{ass}}

\newcommand{\RN}[1]{\textup{\uppercase\expandafter{\romannumeral#1}}}

\newcommand{\rim}{\textcolor[rgb]{0.00,0.00,1.00}}
\def\L{b}
\newcommand{\R}{\mathbb{R}}

\def\How{However, }
\def\bep{\begin{pmatrix}} \def\eep{\end{pmatrix}}
\newcommand\CRNs{chemical reaction networks}

\def\l{{\lambda}}
\def\com{compartment}

\providecommand{\pp}[1]{\left[#1\right]} 
\providecommand{\pr}[1]{\left(#1\right)} 

\usepackage{listings}
\usepackage{color}

\definecolor{dkgreen}{rgb}{0,0.6,0}
\definecolor{gray}{rgb}{0.5,0.5,0.5}
\definecolor{mauve}{rgb}{0.58,0,0.82}

\renewcommand{\theta}{\vartheta}

\renewcommand{\thefootnote}{\fnsymbol{footnote}}

\numberwithin{equation}{section}

\def\satd{satisfied}

\def\bc{\begin{cases}
  }     
\def\ec{\end{cases}}

  \newcommand{\beq}{\begin{eqnarray}
    }
    \usepackage{hyperref}
    \usepackage{cleveref}
\def\eeq{\end{eqnarray}}
   \newcommand{\be}[1]{\begin{equation}\label{#1}}
\newcommand{\ee}{\end{equation}}
\def\bea{\begin{eqnarray*}}\def\ssec{\subsection}
    
\def\eea{\end{eqnarray*}} \def\la{\label}   \def\ith{it holds that } 
      \def\saty{satisfy}        

\def\I{\infty} \def\Eq{\Leftrightarrow}
  \def\T{\widetilde}
  
\def\BEN{\begin{enumerate}}  \def\BI{\begin{itemize}}
\def\EEN{\end{enumerate}}   \def\EI{\end{itemize}} \def\im{\item} \def\Lra{\Longrightarrow}  \def\eqr{\eqref}

\def\mR{\mathcal R} 
\def\sssec{\subsubsection}

\def\g{\gamma}   \def\de{\delta}  \def\b{\beta}
   
\def\ep{\epsilon} 
  
\def\Fr{Furthermore, }

   \def\wrt{with respect to }
  \def\mbw{may be written as }\def\resp{respectively}  \def\fno{from now on} 
 \def\eqr{\eqref}  

\def\a{\alpha}
\def\l{\; \mathsf l}
\def\fr{\frac} \def\im{\item}
\newcommand{\s}{\mathsf s}
\renewcommand{\i}{\mathsf i}
\renewcommand{\r}{\mathsf r}

\def\eeD{\end{definition}} \def\beD{\begin{definition}}
\def\beR{\begin{remark}} \def\eeR{\end{remark}}
\def\beL{\begin{lemma}} \def\eeL{\end{lemma}}

    \def\beP{\begin{proposition}} \def\eeP{\end{proposition}}
    \def\beC{\begin{claim}} \def\eeC{\end{claim}}

    \def\m0{{\mathcal R}_0} \def\w{\omega} 

    \long\def\symbolfootnote[#1]#2{
\begingroup
\def\thefootnote{\fnsymbol{footnote}}\footnote[#1]{#2}
\endgroup}
\def\fn{\symbolfootnote}

\def\l{\lambda}

 \def\D{\Delta}

\newcommand{\red}{\textcolor[rgb]{1.00,0.00,0.00}}
\newcommand{\blue}{\textcolor[rgb]{0.00,0.00,1.00}}

 \def\brn{basic reproduction number } 
\newcommand{\figu}[3]{
\begin{figure}[H]
\centering
\includegraphics[scale=#3]{#1}
\caption{#2\label{f:#1}}
\end{figure}
}

\def\Lgn{\g+\mu+\de} 

\def\brn{basic reproduction number }
 \def\DFE{disease free equilibrium}\def\rd{\r_{dfe}} \def\sd{\s_{dfe}}
 \def\mA{{\mathcal A}}

\begin{document}
\title{Dynamics of an SIR epidemic model with limited
medical resources,  revisited and corrected
}

\author{Rim Adenane$^\text{a}$, Florin Avram$^\text{b}$,  Mohamed El Fatini$^\text{a}$, R.P. Gupta$^\text{c}$}
    \maketitle                         
           \begin{center}
	$^\text{a}$
		Département des Mathématiques,  Universit\'e Ibn-Tofail, K\'enitra, 14000, Morocco,\\
		$^\text{b}$
		Laboratoire de Math\'{e}matiques Appliqu\'{e}es, Universit\'{e} de Pau, 64000, Pau,
 France,\\
  $^\text{b}$ 
  Institute of Science, Banaras Hindu University,Varanasi, 221 005, India
	\end{center}

\begin{abstract}
This paper generalizes and corrects  a famous paper (more than 200 citations) concerning Hopf  and Bogdanov-Takens bifurcations due to L. Zhou and  M. Fan, ``Dynamics of an SIR epidemic model with limited medical
resources revisited",   in which we discovered a significant numerical error.
Importantly, unlike the  paper of Zhou and   Fan and several other papers which followed them, we  offer  a  notebook where the reader may  recover   all the results,  and also modify  them for analyzing similar models. Our calculations lead to the introduction of some interesting symbolic objects, ''Groebner eliminated traces and determinants" -- see \eqr{dt}, \eqr{Gdt}, which seem  to have appeared here for the first time, and which might be of independent interest.
We  hope our paper might serve as yet another  alarm bell regarding the importance of accompanying  papers involving complicated hand computations by  electronic notebooks.
\end{abstract}



\section{Setting the problem}

 This work is motivated by the  task  of obtaining stability and bifurcation results
for a  ``Capasso-Ruan-Wang"  SIR-type epidemic model
\be{Vys}
\bc
\s'(t)=
\pp{\L -\mu \; \s(t)} -\s(t) \pp{\g_s    + \i(t) N(\i(t))}
+i_s \i(t)+ \g_r \r(t),
\\
\i'(t)= \i(t) \pp{
 \s(t) N(\i(t)) -(\g +\mu_i)}-T(\i(t)), \quad \g=i_s+i_r, \mu_i=\mu + \de\\
\r'(t)= \g_s \s(t) + i_r   \i(t)  -{(\g_r+\mu)} \r(t) +T(\i(t)),
\ec
\ee 
which
generalizes many  models studied in the last decades, including works by Capasso and Serio,
 Shigui Ruan, Wendy Wang, L. Zhou and  M. Fan, Vyska and Gilligan, E. Avila-Vales, Pei Yu and coauthors \cite{Capasso,WR,Wang,ZhouFan,Vyska,Rivero,Perez,LuRuan21,Pan21,Gupta}. 

 This model includes
\begin{itemize}
\im {Known  rates:}  $ \L,\mu,\de,\g_s \;$ (of births, deaths, deaths due to epidemics, and vaccination).
\im  {Three statistically estimable transition rates:} $\g_r,  i_r, i_s  $
 (of loss of immunity, and of recovery rates with/without immunity).

\im A nonlinear correction $N(\i)$ to the usual  constant infectivity rate, introduced and studied in the influential papers by Capasso \cite{Capasso},  by Liu et.al. \cite{Liu87}, and by Hethcote and Van den Driessche \cite{hethcote1991some}, with the purpose of capturing the {\bf psychological reaction} of a population to the evolution of the epidemics (note that this type of terms is also popular in the ecology literature \cite{Wol}).

\im A nonlinear, bounded treatment  $\R_+ \ni \i \rightarrow T(\i)\in \R_+$, which unifies
 two cases previously studied in the literature:
 $$T(\i)=\bc\eta \i \quad 0\leq \i\leq \w\\
 \eta \w \quad \i>\w\ec \; , \quad T(i)=\eta \fr{\i}{1+  \i/\w},\quad  \eta>0, \w>0.$$
  Note that both parametrizations have been chosen to  behave as $$T(\i)\approx \bc \eta {\i}&\i\approx 0\\\eta \w&\i\approx \I \ec.$$

Here, we only consider the second case (which may be viewed as a smooth approximation of the more natural first case).
\end{itemize}

\beR
The  vaccination rate of the susceptibles $\g_s$  is absent in the previous works.  
  We also generalize by assuming the  infected individuals may either  transition to the recovery compartement at rate $i_r$, after successful treatment,  or    revert to susceptible  at rate $i_s$, without ever becoming immune (the treatment could be  partially successful, but fail to produce immunity).

\eeR

\beR

  We are using in \eqr{Vys}  a unified notation scheme proposed in \cite{AABH,AAH,AABBGH}, which could be applied to any compartmental model. We propose that
{ a  linear rate of transfer} from compartment $m$ to compartment $c$ be denoted by $m_c$, and the total linear rate out of $m$ is denoted by $\g_m$,  which implies $\sum_{c} m_c= \g_m$ (for example the total removal rate is  $\g_i=i_s+i_r$), and
 new infection parameters be all denoted by $\b_{c_1,c_2}$, where $c_1$ and $c_2$ are the source and destination compartments. However,    simplifications will be made like $\b_{s,i}=\b, \g_i=\g$ for  notations already well established traditionally.

\eeR

\beR  Assuming  $T(\i), N(\i) \ge 0, T(0)=0$, implies for polynomial dynamical systems the ``essential  nonnegativity" of  \eqr{Vys}, i.e. the fact that it never leaves the nonnegative octant invariant -- see \cite{haddad2010, farina2000}.\\
\eeR

\subsection{The purpose of bifurcation analysis}

 According to John Guckenheimer,  Scholarpedia, 2(6):1517, 
the first goal of bifurcation theory is to produce parameter space maps  that divide the  parameter space into regions of topologically equivalent systems.


This is a very challenging process for  models with many parameters, and  a reasonable way to proceed is to decompose it in the following steps:

\BEN \im  Fix    all of the parameters except two,  in a way that produces an interesting two -dimensional partition of the  space of the remaining  parameters into regions with different numbers of fixed points and stability properties. The determination of interesting  regions requires either some symbolic knowledge of the system, or numeric continuation software ( the most well-known open-wares for the latter   are Python, Julia and  MatCont).
\im Identify the ``corner/co-dimension 2 points" (for example,  the  Bogdanov-Takens points), where more than two regions meet, which may be used to organize the computations.
\im Elucidate, via time and phase-plots, the topological behavior of the dynamics in each of the regions encountered, paying special attention to bifurcations, i.e. to the topological changes which occur when crossing boundaries between regions, and notably the  corners. This may be achieved by any CAS, and    we have found  convenient to use  the Mathematica package EcoEvo -- see  \cite{ecoevoG}.
\EEN

As mentioned, there is a vast literature on variations of our class of  models \eqr{Vys}, reviewed briefly in section \ref{s:his}. However, only a small part follow the three-step philosophy mentioned above, and an even smaller part mention the time-length of the computations,  which has long been recognized  as the essential bottleneck in the parallel disciplines of \CRNs, population dynamics, ecology, etc. For some places where this point is touched on, see  \cite{pachter2005algebraic,macauley2020case}, which propose unifying all these disciplines under the name of
{computational/algebraic biology}.

One of the results reflecting the three-step philosophy is  Figure \cite[Fig. 6]{ZhouFan}, reproduced with our notebook as Figure \ref{f:fig6ns} below,  which deals  with the  particular case  $ \g_r=\g_s=i_s=0$ of \eqr{Vys}, and uses the parameter $\alpha$ of \cite{ZhouFan} defined by $\alpha:=\eta \omega$.

\figu{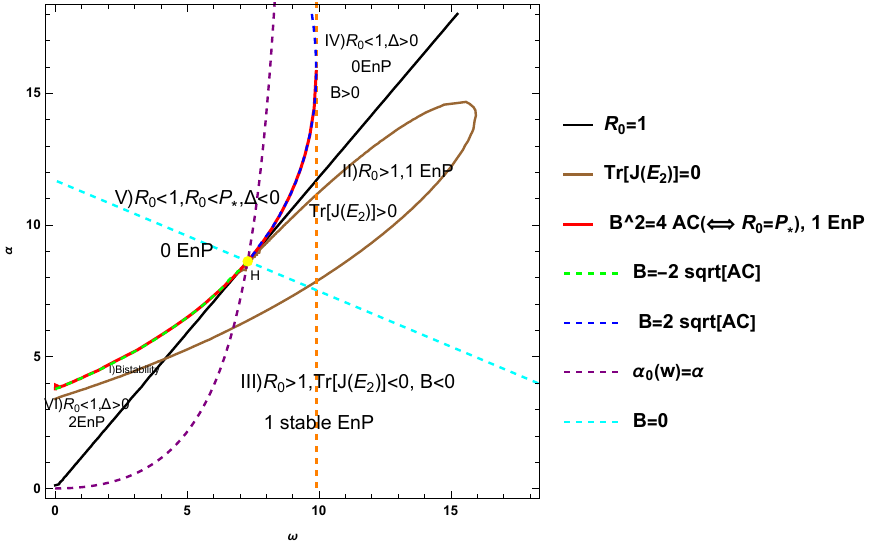}{The original  partition in \cite{ZhouFan} of the  $(\omega,\alpha)$ plane into  regions, defined by the $\m0=1,Discriminant=0, Tr\pr{J_{E_2}}=0$ (other less important dotted curves from the original  are included as well). The   ``yellow point" at which five curves seem to meet turns out to be an optical illusion, which covers an interesting seventh region, invisible at this scale.}{1}

We are building below on the results of  \cite{ZhouFan}. This paper is very well written;  it should be mentioned however  that it is written somewhat from a ``I don't need a computer/high school olympiad problem"  point of view, in the sense that the authors dedicate  pages to switching to new parameters and explicitizing  implicit equations,  
which is typically unnecessary when using a computer.


Our first improvement on \cite{ZhouFan} is the discovery that Figure \cite[Fig. 6]{ZhouFan} contains an optical illusion. Due to the choice of parameters, three important points appear superposed in the yellow point in the figure \ref{f:fig6ns} (denoted by $e$ in \cite{ZhouFan}), where it appears wrongly that five curves meet.

This  is
 corrected in our pictures below \ref{f:fig6ns}, \ref{f:fig6BT}. Note that producing this small modification  translated into  months of work for us, due to the omission by \cite{ZhouFan} to provide their notebooks.

\ssec{The importance of  computational aspects in the study of dynamical systems}
It may be argued that nowadays, all papers written on    applied dynamic models should be  approached from a computational science point of view.
Indeed, the fundamental  problems of determining the  number of fixed points and limit cycles,  and of  identifying  bifurcations (Hopf, saddle-node, Bogdanov-Takens), etc,  reduce all to the solution of polynomial systems of  equalities or inequalities. Now  the modern tools  for  tackling symbolically such problems come  from algebraic geometry (like Groebner bases).

Sometimes, CAS's (which use Groebner bases for solving) solve our problems   immediately, sometimes one needs to   ``twist"  the  notebook before reaching a result, and  sometimes one needs to plug in a numeric condition (possibly with parameters  determined by physical experimentation).
In all these cases, the notebook becomes a crucial part of a paper.

The remarkable  opportunity we have nowadays of being able to accompany our pencil calculations with electronic notebooks was emphasized already 30 years ago in  papers like
``Electronic documents give reproducible research a new meaning" \cite{Repro} -- see also \cite{buckheit1995wavelab,Donoho}. Following the efforts of numerous people,  lots of progress has been achieved in this direction, as witnessed in particular by the existence of the platform GitHub.

\subsection{ A critic of the mathematical epidemiology bifurcations literature}

 We believe that a great part of the mathematical epidemiology literature  suffers today from two   draw-backs.
\BEN \im   The major one  is the lack of ``electronic reproducibility", i.e.  of  notebooks made public.
Here, we provide  electronic notebooks   for the papers cited above,  where the reader may  recover the results,  and also modify  them as he pleases, for analyzing similar models.
We believe that electronic  reproducibility in the above sense  should    be  the norm in the  field of mathematical biology.

 \im A  second draw-back in our opinion is the focus on particular cases, furthermore studied under simplifying assumptions, which are  often  difficult to justify epidemiologically, like assuming $0$ deaths or eternal immunity. It seems sometimes  that  the   justification of the simplifying assumptions is to  allow a ``solution by hand", in the spirit of ``mathematics olympiad problems". Solving particular cases may of course be quite useful pedagogically, but is not efficient.  For just one of many examples where symbolic computing beats easily hand computations,
 note that the ``mathematics olympiad literature" restricts itself often to polynomials  of order $2$, to be able to use their elementary discriminant formula, which is of course unnecessary nowadays. 
 We believe that mathematical epidemiology
 could be better served by relying more heavily on  modern symbolic computation. As another benefit,  from this point of view  there is no need  to dedicate separate papers  to each  particular case  of general models like \eqr{Vys}.  A more natural approach   would be to study unified general cases symbolically as far as possible (until encountering time constraints),   pinpointing  how far can one push symbolics, and  turn to symbolic or  numeric particular cases only subsequently.
 
 Below, we investigate our model at two levels:
 \BEN \im The general one in \eqr{ZFG}, for which we may already provide simple information, like $R_0$, and the maximum possible number of fixed points.
 \im The particular case \eqr{ZF}, for which we are able to provide bifurcation results  as well.

 We believe that the traditional paradigm of one model, one paper doesn't suit well  mathematical epidemiology, where the interesting model for epidemiologists is always quite  general, but mathematicians  are only able to solve completely  just particular cases (which may be uninteresting  for epidemiologists).
\EEN

\EEN

{\bf Contributions}. Our motivation for   revisiting the bifurcation work of  \cite{ZhouFan}  was   the desire to accomplish the  task of
{ providing an electronic notebook, which may eventually be easily modified for studying variations of similar models}.
This quite non-trivial task required in particular pinpointing the ``symbolic-numeric boundary" for this problem, which turned out to be the computation of the trace at the stable endemic point $E_2$. This calculation was first achieved in numeric instances, but finally also symbolically by Mathematica, and lead to the introduction of some interesting symbolic objects, ''Groebner eliminated traces and determinants" -- see \eqr{dt}, \eqr{Gdt}, which seem  to have appeared here for the first time.

{\bf Contents}.
Our paper starts in section \ref{s:his} highlights of  previous  literature concerning the SIRS model with nonlinear force of infection and treatment. 
It then formulates in section \ref{s:ZFG}, \eqr{ZFG}, a general epidemic model, for which $R_0$ and the ``scalarization" of the fixed point system are easily achieved. The analysis of the simplest particular case  $ \g_r=\g_s=i_s=0$, studied in \cite{ZhouFan,Gupta}, is then revisited and corrected in Section \ref{s:ZF}.  
Section \ref{s:Gro} introduces the interesting  ''Groebner eliminated traces and determinants"  \eqr{dt}, \eqr{Gdt}, which turn out useful in the determination of Bogdanov-Takens bifurcations. Sections \ref{s:Fanb}, \ref{s:tpp} complete  aspects missing in the original papers concerning the two-parameter bifurcation diagram and the  important ``corners where several regions meet".
Last but not least, our  paper is accompanied   by electronic notebooks supporting the results.

 \section{A  review of the literature on the SIRS model with nonlinear force of infection and treatment \la{s:his}}

 The epidemic model \eqr{Vys} is inspired firstly by the pioneering work of Capasso and Serio \cite{Capasso} (having currently 1100 citations), which propose to capture in $N(\i)$ psychological (smooth) self regulating adaptations of the population to the epidemics.
Further
 works of Liu, Hethcote,
  Levin and Iwasa \cite{Liu86,Liu87} {on power type interactions $\s \i N(\i)\sim \s \i^p$} suggest adopting a classification of forces of infection   
 in three cases, $p<1,p=1,p>1$. Subsequently, researchers have turned to studying fractional forms like {$\s \i N(\i)=\s \frac{\b \i^{\epsilon}}{1+\xi_1 \i +\xi_2 \i^2} ,\; \epsilon\in \mathbb{N} $}.
  We recommend the choice $2$ for the degree of the denominator, since this  allows, by varying $\ep$, for including all the three cases of \cite{Liu86,Liu87}. The corresponding cases become now: \BEN \im $\ep =1$ (``non-monotone case"),  see for example  \cite{Pan21}, \im $\ep \ge 3$ (``ultimately increasing"), see for example  \cite{LiLi}, and  \im $\ep =2$ -- see \cite{RuanWang,Ruan08}.
  \EEN

In this paper, we will only consider the case $\ep=1$ \fno.\\


 The pioneering papers of Wendy Wang and Shigui Ruan \cite{WR,Wang} initiated the study of  the important  treatment term $T(i)$, which attempts to capture
 a possible partial break-down of the medical system due the epidemics, and \cite{zhang2008backward} proposed to study also    smooth treatment terms like $T(\i)=\eta \fr{\i(t)}{1+ \nu \i(t)}$, motivated by ``the effect of delayed treatment when the number of infected individuals is getting large".
The Capasso-Ruan-Wang model \eqr{Vys} generated lots of works focused on understanding the possible bifurcations which may arise. For a selection of further   works out of a huge literature, see \cite{Wang,Ruan12,ZhouFan,Rivero,Vyska,Jana,Perez,Wei21,LuRuan21,xu2021complex,Pan21,Gupta}, and references cited therein.\\

 Our goal here is to     complete these works  with Mathematica notebooks.

 There are two streams of literature, with smooth and non-smooth (piecewise) treatment, which we aim to  compare. We will start in this first instalment with the smooth case, which is easier to present.

  On a second level, we encounter two levels of complexity in the papers above. The simpler flavor, with $\g_r=0$ or $\de=0$, allows reducing to two dimensions (for example, $\s,\i$), and is encountered in \cite{ZhouFan,Rivero,Wei21,LuRuan21,Pan21,Gupta}. These works, which are all small variations of \eqr{Vys} with $\l=0=\g_r$,   will be reviewed in this first instalment. The more complex three-dimensional version encountered in Z. Hu, W. Ma, S. Ruan (2012) \cite{Ruan12} and Y. Xu, L. Wei, X. Jiang, Z. Zhu (2021) \cite{xu2021complex} will be discussed in a future paper.

We review now briefly the works of \cite{Vyska,Wang,ZhouFan,Jana,Pan21,Perez}.
  In Vyska and Gilligan (2016) \cite{Vyska}, the   infectivity is constant, and
the treatment is piecewise    $$T(\i(t))=\eta  \min[\i(t),\w ],$$ where $\eta=T'(0)$ is  the treatment rate when $\i=0$,
and  $\w $ is an upper limit for the capacity of treatment.

Wendy Wang (2006) \cite{Wang} and Pan,  Huang, and  Huang (2021) \cite{Pan21} consider  piecewise treatment rate, and  diminishing infectivity functions given respectively by
$${N(\i(t))}=\bc \b \fr{1}{1+ \xi \i(t)}\\\b \fr{1}{1+ \xi \i^2(t)}\ec, \xi \ge 0;$$
these reflect a self-regulating  behavior of  the infected people, which reduce their social contacts, in the latter case to $0$ (in the limit $\i \to \I$).  Such negative feedback
effects are called  Holling-type  in ecology, and Michaelis-Menten functions in chemistry and molecular biology.
Note that $\xi=0$ recovers the standard force of infection of the first paper.

Both the above functions, as well as their unification ${N(\s(t)\i(t))}=\b \fr{1}{1+ \xi \i(t)^{\ep}}, \xi \ge 0, \ep >0,$
  satisfy the conditions (I)--(V) in \cite[p.(47)]{Capasso}.

The papers \cite{ZhouFan,Rivero,Jana,Wei21,Gupta}, our focus below, are all particular cases of \eqr{Vys} with $\l=0=\g_r$ (the smooth treatment term is sometimes written as
$T(\i)=\a \fr{\i(t)}{1+ \nu \i(t)}$,
where $\nu= \fr 1 \w, \a =\eta \w$).

 The paper \cite{Rivero} assumed  $i_s>0$, following \cite{Ruan12}, who had introduced this parameter for the non-smooth treatment case. Both these papers work within  the rather restrictive constant population case, by assuming $\mu=\de=0$.

  The  paper \cite{Wei21} considers also a stochastic model, including both Gaussian white noise and a Markov jump process \cite[(6-7)]{Wei21}. They also reproduce in \cite[Fig. 1]{Wei21} the ``Zhou-Fan" map, in a different and maybe more natural parametrization $(\a,\b)$.
   The paper \cite{Jana}  provides  also control results
(supposing $T(\i)=\eta u \fr{\i(t)}{1+ \nu u\i(t)}, u\in [0,1]$). Finally, the  paper \cite{Perez} extends the previous works by considering logistic growth of the susceptibles.

 \beR

 We draw now attention to a  specific problem raised by the use of extra parameters in
\cite{Ruan12,Rivero}, who denote the single parameter $i_s$ by the product $p \de$, using thus two letters for the same parameter  (they also  use $\de(1-q)$, where $q=1-p$, for $i_s$).
Thus,  their starting model contains unnecessary parameters, which of course will not appear in any of the results (and all the equations would be  shorter to write without these $p,q$).  \Fr they denote $\g_r,\g_s$ by $b m, b (1-m)$.
Let us  emphasize   that  the  bifurcation results obtained in  these  papers  are quite interesting.
 \How without adopting a unified notation like $m_c$
for a transfer parameter from $m$ to $c$,   the growing epidemic literature  will become very hard to evaluate.
The problem is further compounded in  \cite{xu2021complex}, who denote  $\g_r$ by $b m +\b$ (without any concrete attempt to relate this three parameters to statistical data, which might have justified the extra parameters.

 \eeR

\section{Generalization of the \cite{ZhouFan,Gupta} model 
} \la{s:ZFG}
Some of the results of \cite{ZhouFan,Gupta} hold also for the more general model:
\be{ZFG}
\bc
\s'(t)=
\L- \s(t) \pr{\g_s+  \mu +\b \fr{\i(t)}{1+ \xi \i(t)}}+i_s \i(t)+\g_r \r(t)
\\
\i'(t)=
\i(t)\pp{ \s(t)  \fr{\b}{1+ \xi \i(t)} -\eta \fr{\w }{\w+  \i(t)} -v_i} \\
\r'(t)=\g_s \s(t)+ i_r   \i(t)+\eta \fr{\w \i(t)}{\w+  \i(t)}  -(\mu+\g_r)  \r(t)
\ec,
\ee
where $v_i= \Lgn.$ It has  $12$ independent parameters $(\L,\mu,\b,\g,\g_r,\g_s,i_s,i_s,\delta,\xi,\eta,w)$.

\ssec{Preliminary steps}
The qualitative analysis of the model \eqr{ZFG} starts by identifying the ``infectious \com"  $\i$ and its ``non-infectious " complement,  $\s,\r$.
Subsequently,
\BEN
\im The unique ``\DFE" (DFE) boundary fixed point
  \[E_0=(\frac{\L   \left(\mu +\gamma _r\right)}{\mu  \left(\mu +\gamma _r+\gamma _s\right)},0,\frac{\L   \gamma _s}{\mu  \left(\mu +\gamma _r+\gamma _s\right)}):=(\sd,0,\rd)\] is obtained by plugging $\i=0$ in the ``non-infectious equations".

\im
The Jacobian at the DFE  has block form
\bea jac(DFE)=\left(
\begin{array}{ccc}
 -\mu -\gamma _s & i_s-\frac{\beta  \L   \left(\mu +\gamma _r\right)}{\mu  \left(\mu +\gamma _r+\gamma _s\right)} & \gamma _r \\
 0 & -\eta -\g - \mu -\delta +\frac{\beta  \L   \left(\mu +\gamma _r\right)}{\mu  \left(\mu +\gamma _r+\gamma _s\right)} & 0 \\
 \gamma _s & \eta +i_r & -\mu -\gamma _r \\
\end{array}
\right),
\eea
 and the local stability region of the DFE may be written as $\m0=\sd \frac{ \b}{\Lgn+\eta}<1$, where $\m0$ is the famous ``\brn". For the case without treatment $\eta=0$, this reduces to  $\T {\m0}=\sd \frac{ \b}{\Lgn}$.
 At $\m0=1 \Eq \eta=\eta_0:=\fr{\b \L}{\mu} -v_i=v_i (\T{\m0}-1)$ we have the usual trans-critical bifurcation \cite[Thm. 2]{Gupta}.

\im There may be at most two endemic points, whose algebraic expressions may be quickly found by solving the system symbolically. An interesting alternative  is to reduce the system, with the second equation divided by $i$,
to a scalar equation in $i$, and looking subsequently for positive solutions. This is  standard procedure in the field, and there are in fact cases where this procedure works quicker than the global ``Solve", due probably to the help we are giving the CAS by indicating what to eliminate. The natural strategy for elimination is to use Groebner bases. Our utility  Grobpol reduces a system  to a  scalar polynomial in ``ind", a variable to be specified by the user:
\begin{verbatim}
 Grobpol[mod_, ind_, cn_ : {}] := Module[{dyn, X, par, eq, elim},
         dyn = mod[[1]]; X = mod[[2]]; par = mod[[3]]; 
         eq = Thread[dyn == 0]; elim = Complement[Range[Length[X]], ind];
         pol = 
         Collect[GroebnerBasis[Numerator[Together[dyn /. cn]], 
         Join[par, X[[ind]]], X[[elim]]], X[[ind]]];
         ratsub = Solve[Drop[eq, ind], {s, r}][[1]]; {ratsub, pol}
      ]
 \end{verbatim}

 This  reduces our system to a fourth order polynomial which factors into $i \times$, a linear term with negative root, and a quadratic polyomial in  $\i$
\be{qi} p(i)= A i^2+B i + C=0,\ee
with coefficients  computed at the end of the first cell in \cite{NotM1}. They recover in particular the coefficients of
 \cite[(2.4)]{ZhouFan}:
 \be{ABCF}
 \bc
 A = v_s  v_i >0,\\
 B=    V_i  v_s  \w  -\mu \eta_0  \\
  C= \w \mu  V_i \pr{1-\m0},
 \ec
\ee
 where we put \be{v}\bc
 v_s=\beta +  \xi \mu\\
 V_i=v_i+\eta. \ec\ee

We will order the two possible endemic points $E_1,E_2$ by their $i$ coordinates, given by
 $ i_1= \frac{-B-\sqrt{\Delta}}{2A}  \; \mbox{and} \quad i_2= \frac{-B+\sqrt{\Delta}}{2A},$  where
 \be{dis} \Delta=B^2 -4 AC= \alpha ^2 v_s^2+2 \alpha  v_s (\omega  v_s v_i-\mu  v_i-\beta  \L  )+(\omega  v_s
  v_i-\mu  v_i+\beta  \L  )^2\ee denotes  the discriminant of  the equation. Note that when $C = 0\Eq \m0=1$ we have $ i_1=0, i_2= -\frac{B}{A}$, which indicates a transition in the number of  endemic points when crossing this boundary, and that the equation $\Delta=0$ \mbw\
  $$\pp{\alpha  v_s + (\omega  v_s v_i-\mu  v_i-\beta  \L  )}^2=4 \b \L v_i(\mu-\w v_s)$$
 has solutions involving square roots, with respect to any of its parameters.

\EEN

\section{Review and corrections of the results of \cite{ZhouFan}
}\la{s:ZF}

We  revisit now \cite{ZhouFan,Gupta}, where $\g_r=\g_s=i_s=0$.
The system \eqr{ZFG} becomes
\be{ZF}
\bc
\s'(t)=
\L- \s(t) \pr{  \mu +\b \fr{\i(t)}{1+ \xi \i(t)}}
\\
\i'(t)=
\i(t)\pp{ \s(t)  \fr{\b}{1+ \xi \i(t)} -\eta \fr{\w }{\w+  \i(t)} -v_i} \\
\r'(t)= \g   \i(t)+\eta \fr{\w \i(t)}{\w+  \i(t)}  -\mu \r(t)
\ec,
\ee
where $v_i= \Lgn.$ It has therefore $8$ independent parameters $(\L,\mu,\b,\g,\xi,\eta,w,\delta)$.

Since  $\g_r=0 $ implies that $ \r$ doesn't appear in the first two equations in \eqr{ZF},  it is enough to   study  the two dimensional system for $(\s,\i)$, with Jacobian

\be{J}
J= \left(
\begin{array}{cc}
 -\frac{\beta  i}{i \xi +1}-\mu  & -\frac{\beta  s}{(i \xi +1)^2} \\
 \frac{\beta  i}{1+i \xi} &  -\frac{\eta  \omega ^2}{(i+\omega )^2}+\frac{\beta  s}{(i \xi +1)^2}-v_i  \\
\end{array}
\right)=\left(
\begin{array}{cc}
 -\frac{\L}{s}  & -\frac{\beta  s}{(i \xi +1)^2} \\
 \frac{\beta  i}{i \xi +1} &  \frac{\eta  \omega i}{(i+\omega )^2}-\frac{\beta \xi s i}{(1+ \xi i)^2}  \\
\end{array}
\right),
\ee
where the second form, used in \cite[Lem. 2]{Gupta},  
holds only at the endemic points.

\ssec {The stability of the endemic points}

   Establishing  the stability of the two endemic points in \cite[Thm 3.2-3.3]{ZhouFan} is rather tricky, since it requires working with the complicated quantities
 $det(J(E_1)),  det(J(E_2)), $
 $Tr(J(E_1)), Tr(J(E_2))$. 
 One possible approach to efficient symbolic computations in this example is keeping the variable  $i$, but eliminating $s$ either  from the first stationarity equation in \eqr{ZF}, yielding $s=\fr{\L(1+\xi i)}{\mu + v_s i}, v_s=\b+\mu \xi,$ or from the second equation in \eqr{ZF}, which yields $s=\frac{(i \xi +1) \left(\eta  \omega +i v_i+v_i \omega \right)}{\beta  (i+\omega )}$.

Using the second  equation, one finds that at the endemic points
\bea det(J)= \fr{i}{(1+i \xi)(i+\omega )^2} \Psi(i), \; \; \Psi(i)=v_s v_i (i+\omega )^2+\a  (  \omega v_s-\mu ), \eea which checks with the   second order polynomial which appears  in the proofs of \cite[Thm 3.2-3.3]{ZhouFan}. The trace after eliminating $s$ from the second stationarity equation also confirms the result on \cite[pg. 319]{ZhouFan}.

  {Cf. \cite[Lem. 2,Thm 3]{Gupta}}
   \be{3d} det(J(E_1)) < 0 < det(J(E_2)),\ee
which implies that  $E_1$ is a saddle whenever it exists \cite[Thm 3.2]{ZhouFan}, and that the stability region for $E_2$
 coincides with the region where $Tr(J(E_2))<0$ \cite[Thm 3.3]{ZhouFan}, \cite[Thm. 3]{Gupta}.

Now  \eqr{3d} is immediate when $\w v_s >\mu$, because in that case $\Psi(i)>0,i_1 <0,$ and $i_2>0$, whenever it exists.   For the other case, see the elegant proof in \cite[Lem. 2,Thm 3]{Gupta}.

  Since $E_1$ is a saddle and the determinant at $E_2$ has constant sign, it follows that only three hyper-surfaces  are necessary for  bifurcation analysis:
\be{BTZF}\bc \m0=1 
\\\Delta=B^2 -4 AC=0
 \\Tr(J(E_2))=0\ec.\ee

 The next task is determining  symbolically, or, at least numerically,  the intersection of these curves.

As mentioned already,  the expression of the trace $Tr(J(E_2))$ after explicitizing  $i_2$ is very complicated. At this point, one has two possibilities to proceed with:
 \BEN  \im A symbolic attack on the trace $trG$, without choosing a fixed point by explicitizing $i$, via a  GroebnerBasis command

\begin{lstlisting}
// trG=GroebnerBasis[{dyn,Tr[jac]},Pars,{s,i}]
\end{lstlisting}
 which eliminates the variables $\{s,i\}$ from the Trace of the Jacobian.
 This allows producing  a plot of both branches of the  curve   $\pp{Tr\pr{J_{E_2}}=0}\cup \pp{Tr\pr{J_{E_1}}=0}$, together ( see Figure \ref{f:fig62}) (and  includes our desired curve ${Tr\pr{J_{E_2}}=0}$).  This approach is explained in more detail in Section \ref{s:Gro}.

 A  similar alternative approach is  to compute the functions $ tr(i), p(i)$ which depend on $i$ but not on $s$, and their resultant.  Both these approaches work very quickly, but do not allow separating the two branches. This may only be achieved  via the  numeric approach described next.

 \im Explicitize $i=i_2$ and continue with numeric values for  all parameters,  except the two chosen to be displayed  (using rational numbers, to be able to work with infinite precision) -- see Figure \ref{f:fig6ns}.

\EEN

 We report on both approaches in next section, as well as in  a ``symbolic notebook" containing whatever could be computed symbolically: fixed points, bifurcation varieties, and their intersections. This produces the two dimensional map obtained by evaluating our curves under a ``cut  numeric condition" for the remaining parameters,  and examines examples of points in most of the regions, and on their boundaries. Here we must clarify that a very
 important part of a numeric bifurcation project involving many parameters is to obtain an interesting cut where all the parameters but two are fixed. This may be achieved
by packages for ``numerical continuation and bifurcation" of ODE dynamical systems like  MatCont (written in Matlab), PyDSTool (Python),  XPPAuto (C), and BifurcationsKit (Julia)  -- see \cite{blyth2020tutorial} for a recent review.
In our paper, we just perturbed a bit the cut furnished in \cite{ZhouFan}, 
until it became clearer that the coincidence of the four curves in their paper was an optical illusion.

 Finally, we offer  a second ``numeric notebook" which uses the Mathematica contributed  package  EcoEvo, which helps to produce quickly  the time and phase-plots illustrating the various points chosen for display.

 \subsection{A key tool: computing Bogdanov-Takens bifurcations using symbolic algebra}\la{s:Gro}

 Symbolic stability and bifurcation analysis are very time consuming, since they require identifying  varieties like traces, determinants and Hurwitz determinants, evaluated at all the fixed points. It is plausible however that the equation satisfied by the union over the fixed points (i.e. the product of all the  respective equations), ends up simpler symbolically.

  Before presenting the bifurcation results for our example, we present  now  some useful symbolic objects which we have introduced, and were not able to find  in the previous dynamical systems literature.

  \beD \la{d:E} The  determinant, trace, and Hurwitz determinant  of an algebraic system with respect to a subset of its solutions $A$ is defined by  the expressions

 \be{dt} \bc detE=\prod_{i \in \mA} Det(J(E_i))\\ trE=\prod_{i \in \mA} Tr(J(E_i))\\H_{n,E}=\prod_{i \in A} H_n(E_i).\ec\ee
 \eeD

  \beR \BEN \im Two choices of $\mA$ are of special interest: one, involving the product of the traces over all the fixed points, is useful for general dynamical systems. Another one, involving only the product  over  the interior points, is the one we used mostly, since eliminating the DFE from the study is advantageous computationally, and anyway the DFE is easy to study separately.

  \im  The varieties corresponding to $detE,trE$, etc, are the union of the varieties for each branch.  Due to the Vieta relations between the roots, they are expected to be considerably simpler than the individual expressions for the various branches.

  \im Interestingly, $detE,trE$, etc seem related up to proportionality  to the results obtained  via the  GroebnerBasis elimination command

 \be{Gdt}\bc {detG}=\text{GroebnerBasis}[Numerator[Together[[\{\text{dyn}[[1]],\text{dyn}[[2]], Det(X) \}]],\text{par},X]\\{trG}=\text{GroebnerBasis}[Numerator[Together[[\{\text{dyn}[[1]],\text{dyn}[[2]], Tr(X) \}]],\text{par},X]\ec,\ee
 which are easily obtained via the script:
\begin{verbatim}
 GBH[mod_,scal_,cn_ : {}] :=
     GroebnerBasis[Numerator[Together[[Join[ mod[[1]], {scal}]]],
     mod[[3]],  mod[[2]]]
\end{verbatim}
(recall that mod=\{dyn,var,par\}).
  To resolve a two dimensional model,  it is enough to apply this script  with scal= trace and scal=determinant.

 \im The proportionality might be caused by {spurious  implementation factors.}  It may be interesting to clarify the relation between ${detE}, {trE}$ and the corresponding Mathematica objects ${detG}, {trG}$.
\EEN

\eeR

\beR \la{r:easy} \BEN \im The computation of candidates for the Bogdanov-Takens bifurcations for all the fixed points in two dimensions
 may be achieved a priori via the single command\\
 {
 \be{Easy} Reduce[\{trG == 0, detG == 0,cp\}],\ee
 where cp are positivity and other eventual constraints on the variables and parameters.
}

For the \cite{ZhouFan} model, we find in the third cell of \cite{NotM} that detG is proportional to the product of $R_0-1$ and the discriminant, and that trG is proportional to the product of 
\be{tr3} \bc \omega -\fr{(\gamma +\delta +\mu ) \left(2 \mu  (\gamma +\delta +\eta )+(\gamma +\delta ) (\gamma +\delta +\eta )+\mu ^2\right)}{ (\gamma +\delta +\eta +\mu ) \left(\beta  \eta +\xi  \left(2 \mu  (\gamma +\delta +\eta )+(\gamma +\delta ) (\gamma +\delta +\eta )+\mu ^2\right)\right)},\\\xi  \omega ^2 (\beta +\mu  \xi )^2 (\gamma +\delta +\eta +\mu )-\omega  (\beta +\mu  \xi ) (\beta  \eta +\mu  \xi  (2 \gamma +2 \delta +\eta +2 \mu ))+\mu ^2 \xi  (\gamma +\delta +\mu )
\ec \ee  and  a third order polynomial in $\omega$ which has one rational root.

 Since \eqr{Easy} takes too much time, we decompose into $4$ equalities
corresponding to all combinations of the factors.

\im In more dimensions, it is enough a priori to replace the trace in \eqr{Easy} by the maximal
 dimension Hurwitz determinant, and investigate also the positivity conditions associated to the other Hurwitz determinants.

 \EEN
 \eeR

For the generalization of \cite{ZhouFan}, the detG can be obtained, but trG exceeds our computing time (more than one night).  We turn therefore from now on to the particular case of \cite{ZhouFan}.

\subsection{The corrected two parameter bifurcation diagram
} \la{s:Fanb}

We  correct   here the  two-parameter bifurcation diagram in \cite[Fig. 6]{ZhouFan}.
Our  model has  $8$ parameters, or rather seven, since $\de,\g$ intervene in the first two equations only via their sum.
Recall that \cite[Fig. 5-6]{ZhouFan} offer values for   six of the parameters   ($\L= 16,\mu =\frac{1}{10},\delta = \frac{2}{10},\gamma = \frac{12}{100}, (\Lra v_i=\frac{42}{100}), \beta = \frac{1}{100},\xi= \frac{1}{1000} $)  which produced an interesting partition of the two parameter space of the remaining ``treatment parameters" $(\w,\eta)$.
 We have changed however one parameter
to $\mu = .12$ -- see Figure \ref{f:fig6ns}, and \cite{NotM1}, since this makes evident that the point $e$
in \cite[Fig. 6]{ZhouFan}, through which four curves seem to pass, is an optical illusion (which hides a  region VII which is revealed after a further blow-up in Figure \ref{f:fig6BT}).

For comparison with the original, we will show our figure with $\eta$
replaced by the \cite{ZhouFan} parameter $ \a= \eta \w$.

The parameter space is divided in 7 regions, defined by the signs of $\Delta,\m0-1,Tr(J(E_2))$, called \resp\ I,II,III,IV,V,VI,VIa;  the corresponding sign patterns, including that of $ { B}$ (which turns out helpful) are given in the following table:

{\small  \begin{table}[H]
  \begin{center}
    \caption{Table of the sign patterns corresponding to the six regions.}
    \label{tab:table1}
    \begin{tabular}{|l|c|r|r|r|}
    \hline
       & sign of $\Delta$ & sign of $\m0-1$& sign of $Tr(J(E_2))$ & sign of  B \\
  \hline
     Region I ($Q_I$) &  $>0$ & $<0$ &$ >0$ &$ <0$ \\
      \hline
       Region II ($Q_{II}$)&  $>0$ & $>0$ & $>0$& both  \\
      \hline
       Region III ($Q_{III}$) & $>0$  & $>0$  & $<0$& both \\
      \hline
       Region IV ($Q_{IV}$)&  $>0$ & $<0$ & $>0$ & $>0$  \\
      \hline
       Region V ($Q_V$)& $<0$  & $<0$  & $>0$ & both \\
       \hline
       Region VI ($Q_{VI}$)& $>0$  & $<0$ &$<0$ & $<0$ \\
        \hline
        Region VI ($Q_{VIa}$)& $>0$  & $<0$ &$<0$ & $<0$ \\
        \hline
      H& =0& =0&$>0$ & =0\\ 
      \hline
      BT& =0& $<0$& =0&$<0$ \\ 
      \hline
     $B_1,B_2$&$>0$ & =0& =0&$<0$\\ 
      \hline
    \end{tabular}
  \end{center}
\end{table}
}

 \figu{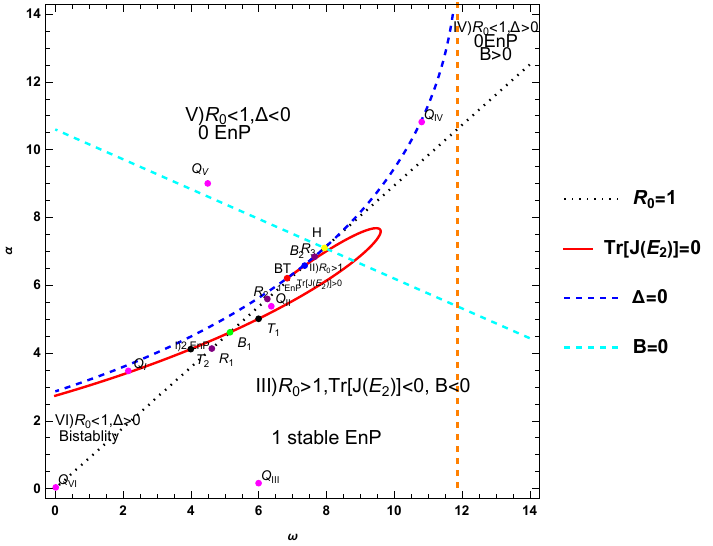}{The partition of the  $(w,\alpha)$ plane into seven regions, defined by $\m0=1,\Delta=0, Tr\pr{J_{E_2}}=0$ (which is involved in determining the Hopf bifurcation variety), {and $B=0$ (the crossing of which changes the number of endemic points  between the regions I,V from $2$ to $0$)}.  Note that the region $ Tr\pr{J_{E_2}}<0$ is included in the region $\Delta \ge 0$, confirming the fact that this yields the stability region for $E_2$
  \cite[Thm 3.3]{ZhouFan}, \cite[Thm. 3]{Gupta}.
  {The trace curve} $ Tr\pr{J_{E_2}}=0$ of possible Hopf points ends at the point $BT$, due to the fact that we broke here {it's twin part} $ Tr\pr{J_{E_1}}=0$ --see Figure \ref{f:fig62}, where no Hopf bifurcations may appear. This curve may further be divided in an upper part separating regions II and III and a lower part separating regions I and VI. When crossing the upper part, the stability of the point $E_2$ and hence of the cycle surrounding it changes. When crossing the lower part from region III into region VI, it was conjectured but not proved in \cite[Fig.6]{ZhouFan}
   that periodic non-attracting orbits  appear in a neighborhood of the boundary. We verify this numerically below in figures \ref{f:IItoIII}, \ref{f:III} and \ref{f:R1}. 
  }{1}

\figu{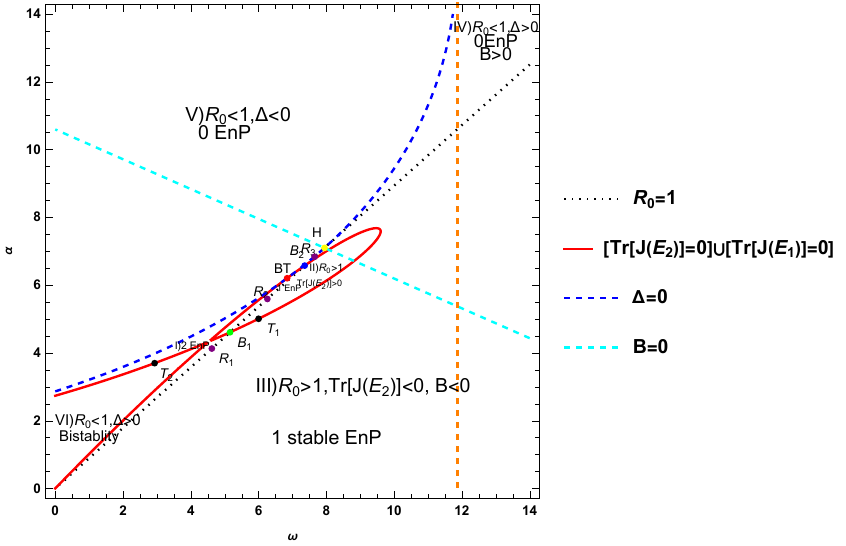}{Same figure, with the curve $ Tr\pr{J_{E_2}}=0$ replaced by the easier to compute $ \pp{Tr\pr{J_{E_2}}=0}\cup \pp{Tr\pr{J_{E_1}}=0}$}{1}

\beR  The point H when both the endemic points collide with the DFE is at the intersection of the  curves $\bc B= 0\\C=0\Eq \m0=1\ec$, where $B,C$ are defined in \eqr{ABCF}. It is easy to determine:
 $ \bc \m0=1 \Lra V_i =\frac{\Lambda}{ \mu} \b \\
  B=0\Lra V_i= \frac{\mu \eta }{v_1 \w}\ec $ yields

  \be{H} \bc \eta_H=\eta_0=\frac{\beta  \Lambda }{\mu }-v_i,\\
  \w_H=\frac{\mu^2 \eta_0}{\beta  \Lambda  v_1} \Lra \alpha_H={\eta_0 \w_H=\frac{\mu^2 \eta_0^2}{\beta  \Lambda  v_1} }  \ec \ee
  The fact that $Tr(E_2)=0$  in Figure \ref{f:fig6F} of \cite{ZhouFan} was an illusion; in fact,  $ Tr_H(E_2)=-\mu$.  However, $\Delta_H=0$, since
  $$\Delta_{\m0=1}=\frac{\left(\beta  \Lambda  \omega v_1 -\mu ^2(\beta  \Lambda/  \mu -v_i )\right)^2}{\mu ^2}=\fr{\beta  \Lambda  v_1}{\mu ^2}\pr{\w -\w_H}^2.$$

   The points at the intersection of $\m0=1$ 
    with $tr(E_2)=0$ may be determined symbolically, for example by computing the resultant of $p(i)$ and the  numerator of $tr(i)$ under the condition $\eta=\eta_0$, and solving it \wrt\ $\w$.\fn[4]{An alternative approach is to   compute the  numerator of $tr(i)$ modulo $p(i)$,  and then determine the sign of the first order remainder evaluated at $i_2$  \cite[Thm. 3.3]{ZhouFan}.}

   Assuming $\T{\m0} \neq 1,\omega \neq \mu /v_1$, this yields that  $\w$ is one of the roots of a {third order polynomial}
   to be called ``B-points".

   In the example illustrated in the figure \ref{f:fig6F} we find, besides the H point,  three B points, and a ``BT" (Bermuda triangle) point, at the intersection of the trace and discriminant curves.
   {We exclude the third B point, since it belongs to the $Tr(E_1)=0$ branch}. The remaining points, together with other points along the boundary $\mR_0=1$,
   ordered by $\w$,  are:
   \bea\bc R_1=\;\{\omega \to 4.620,\alpha \to 4.127\} \mbox{ (on   $\m0=1$, downward from $B_1$)}\\
    B_1=\{\omega \to 5.157,\alpha \to 4.607\}\\
    R_2=\{\omega \to 6.258,\alpha \to 5.590\} \mbox{(on   $\m0=1$, upward from $B_1$)}\;\\
    BT=\{\omega \to 6.841,\alpha \to 6.203\}\\
    R_3=\{\omega \to 7.652,\alpha \to 6.835\} \mbox{ (on   $\m0=1$, between $B_2$ and H)}\;\\
    B_2=\{\omega \to 7.359,\alpha \to 6.574\}\\
    H=\{\omega \to 7.944,\alpha \to 7.097\}\\
    \ec.\eea

    We also display the interior points
    \bea\bc
     Q_I=\{\omega \to 2.156,\alpha \to 3.468\}\\
      Q_{II}=\{\omega \to 6.375,\alpha \to  5.375\}\\
       Q_{III}=\{\omega \to 6,\alpha \to 0.156\}\\
        Q_{IV}=\{\omega \to 11.75,\alpha \to 11.75\}\\
         Q_{V}=\{\omega \to 6,\alpha \to 6\}\\
          Q_{VI}=\{\omega \to 0.078,\alpha \to  0.156\}\ec \eea
          and some  points along the boundary $Tr(E_2)=0$
    \bea\bc
     T_1=\{\omega \to 6 ,\alpha \to 5.00625 \}\\
     T_2=\{\omega \to 2.932 ,\alpha \to 3.696 \}
      \ec \eea

\beR Note the middle two points $BT$ and $B_2$ are hard to distinguish at any scale by eye, since in between them the $\Delta=0,$
its $\m0=1$ tangent at $H$, and $Tr(E_2)=0$  practically coincide. Together with the point $H$, these points give rise to a ``Bermuda triangle, where the phase-plot looks like in Figure \ref{f:fig6BT} below.\eeR

\figu{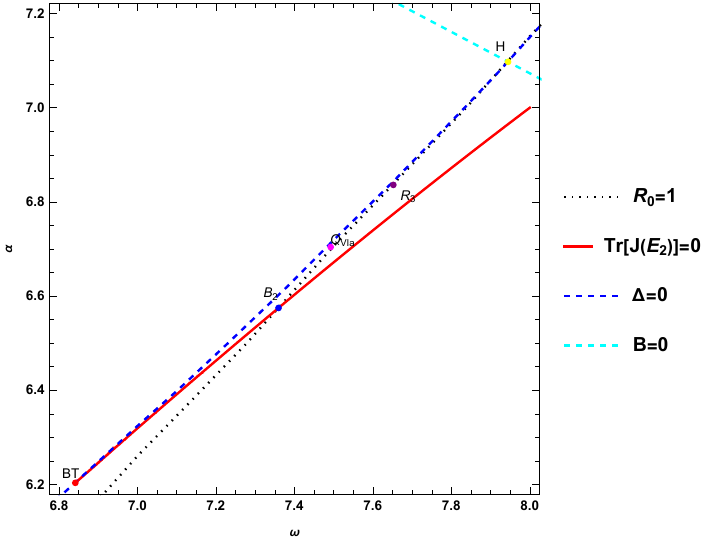}{The Bermuda triangle}{1}

\eeR

\ssec{Some  time and phase plot illustrations of the dynamics in the various regions}\la{s:tpp}

We undertake now a journey between the six regions in Figure \ref{f:fig6ns}, showing  phase and time-plots for specific parameter values.  We switch from now on to using the package EcoEvo in our electronic notebook \cite{NotME}, which has the advantage of computing  also the period of cycles and their Floquet numbers.

We start in the region II, where the {surprising absence of stable
fixed points} (this does not seem to occur in simple models without functional parameters) implies the existence of at least one stable limit cycle.

\sssec{Region II, with all the fixed points  unstable}

A random choice of parameters  confirms that at our first choice (obtained using FindInstance) there exists indeed at least one stable limit cycle, which surrounds the unique  unstable spiral point $E_2=(85.4965, 6.75961)$  with eigenvalues   $0.00887972 \pm 0.13495\; Im  $. The DFE\  $E_0=(133.333,0)$ is a saddle point with eigenvalues  $(-0.12, 0.0501961)$, and the point $E_1$ is outside the positive orthant.

\begin{figure}[H]
    \centering
    \begin{subfigure}[a]{0.4\textwidth}
    \centering
        \includegraphics[width=\textwidth]{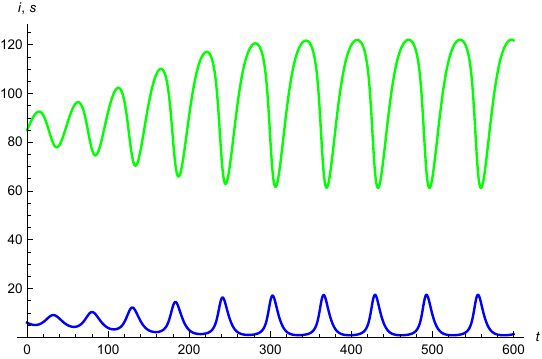}
        \caption{ $(s,i)$-time plot indicate the existence of a limit cycle surrounding $E_2$}
    \end{subfigure}%
    ~
     \begin{subfigure}[a]{0.38\textwidth}
     \centering
        \includegraphics[width=\textwidth]{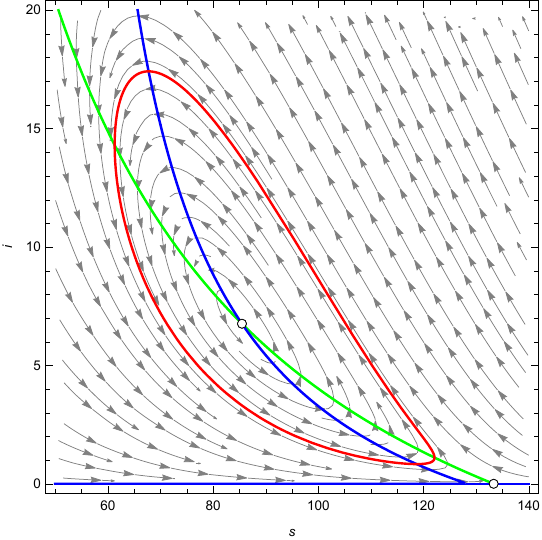}
        \caption{$(\s,\i)$- phase plot.  Attracting orbit with period  $63.5036$, and  Floquet exponents  $\pr{-1.2044\times 10^{-8},-0.0228667}$.}
    \end{subfigure}
    \caption{Time and phase-plot  corresponding to  region II  in Figure \ref{f:fig6ns},with $\left\{ \w=51/8,\alpha=43/8  \right\} \Lra  \m0=1.03912>1 $.}
    \end{figure}

In the next five sections we will investigate the spread of oscillations into the neighboring regions III and I, and at the corner points $B_1, B_2,BTP$.
  \sssec{Spread of oscillations into the neighboring region III}
  
  At the boundary  between regions II and III \ith $Tr(J(E_2))=0$; at $T_1=\left\{ \w=6,\alpha=5.00625  \right\}$,
 $E_2$  is a potential Hopf point. Our first choice, obtained using FindInstance, 
  suggests that the endemic point $E_2=(80.6403,7.90317)$ with eigenvalues  $\pm 0.15125  Im $ is surrounded by  a stable limit cycle.
  The \DFE\  $E_0=(133.333,0)$ is a saddle point with eigenvalues  $(-0.12, 0.0589583)$, and  the fixed point $E_1$ is outside the positive orthant.

\begin{figure}[H]
    \centering
    \begin{subfigure}[a]{0.4\textwidth}
    \centering
        \includegraphics[width=\textwidth]{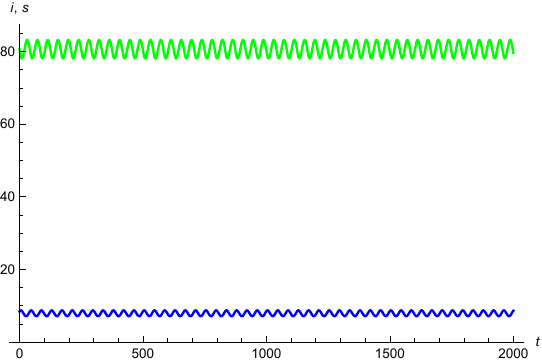}
        \caption{ $(s,i)$-time plot indicate the existence of a limit cycle surrounding $E_2$}
    \end{subfigure}%
    ~
     \begin{subfigure}[a]{0.38\textwidth}
     \centering
        \includegraphics[width=\textwidth]{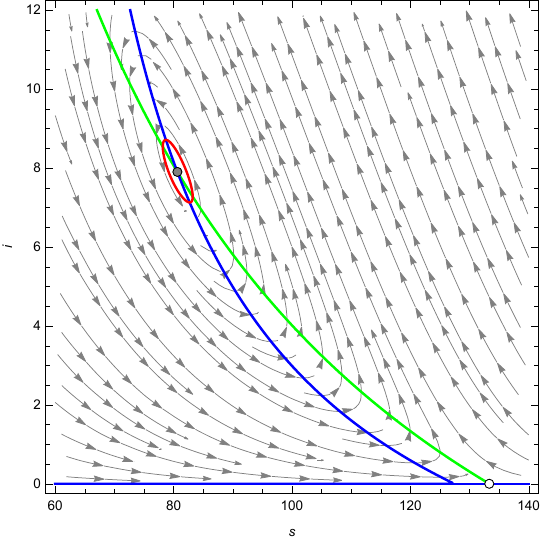}
        \caption{$(\s,\i)$- phase plot.  The period of the cycle is $41.6159$, and the Floquet exponents are $\pr{-3.17375\times 10^{-7} \pm 9.35363\times 10^{-6} \; Im}$}
    \end{subfigure}
    \caption{Time and phase-plot  corresponding to  boundary $Tr(J(E_2))=0$ between regions II and III  in Figure \ref{f:fig6ns}, with $\left\{ \w=6,\alpha=5.00625  \right\} \Lra  \mR_0=1.04626>1$.\label{f:IItoIII}}
    \end{figure}

  Crossing now inside  region III, the  endemic point $E_2$ becomes stable \cite[Thm 3.3]{ZhouFan},
  \fn[4]{It is also globally stable under a certain extra condition \cite[Thm 3.5]{ZhouFan}, while the case when the extra condition is not \satd\ is left open.}

      The cycle is at first unstable "near the boundary", and then disppears ``far enough" from the boundary. The interior point $E_2=(45.58,  23.6495)$  is an attracting  spiral  with eigenvalues   $-0.178557 \pm 0.265984\; Im  $. The \DFE\  $E_0=(133.333,0)$ is a saddle point with eigenvalues  $(0.867292, -0.12)$, and the fixed point $E_1$ is outside the domain.
\begin{figure}[H]
    \centering
    \begin{subfigure}[a]{0.4\textwidth}
    \centering
        \includegraphics[width=\textwidth]{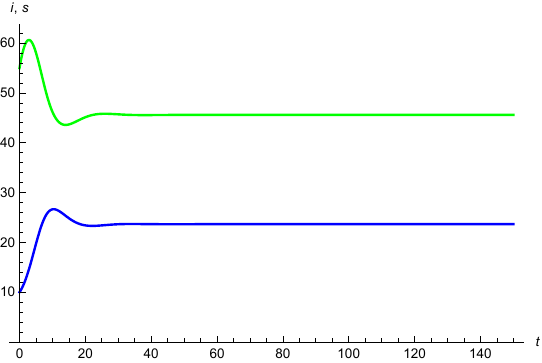}
        \caption{ $(s,i)$-time plot indicate the convergence towards the endemic point $E_2=(45.58,  23.6495)$}
    \end{subfigure}%
    ~
     \begin{subfigure}[a]{0.38\textwidth}
     \centering
        \includegraphics[width=\textwidth]{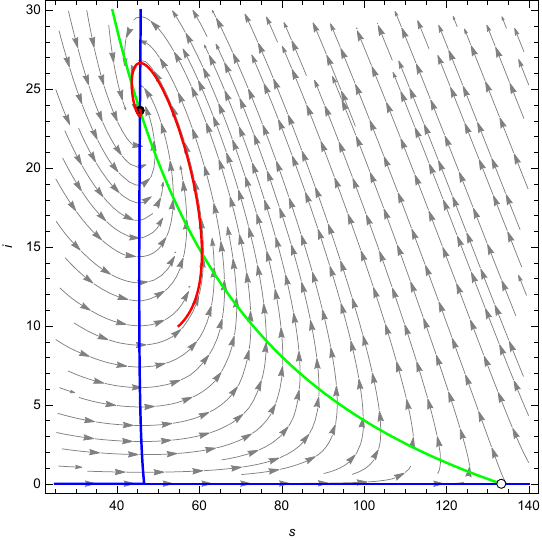}
        \caption{$(\s,\i)$- phase plot. }
    \end{subfigure}
    \caption{Time and phase-plot  corresponding to  region III  in Figure \ref{f:fig6ns},with $\left\{ \w=6,\alpha=5/32 \right\}  $.\label{f:III}}
    \end{figure}

\sssec{Spread of oscillations into the neighboring region I}

Recall that in region II, the absence of any fixed stable points implies the existence of at least one stable limit cycle.
  This is still true on the {boundary $\mR_0=1$ with region I, from $B_1$ to  $B_2$}.

For example, at the point   $\left\{ \w=7.16058,\alpha=6.39679 \right\} $, the following figure reveals the existence of an unstable limit cycle around the endemic point $E_2= (111.34,  2.376) $  with eigenvalues $ 0.0103906 \pm  0.0502167\; Im$. The two fixed points $E_0=E_1=(133.333, 0)$ collide with eigenvalues   $(-0.12,0  )$.
\begin{figure}[H]
    \centering
    \begin{subfigure}[a]{0.4\textwidth}
    \centering
        \includegraphics[width=\textwidth]{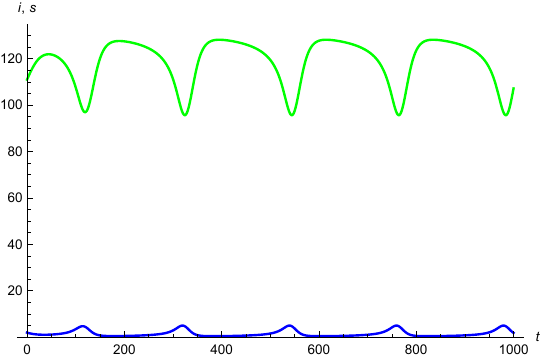}
        \caption{ $(s,i)$-time plot shows the existence of limit cycle surrounding $E_2$ }
    \end{subfigure}%
    ~
     \begin{subfigure}[a]{0.38\textwidth}
     \centering
        \includegraphics[width=\textwidth]{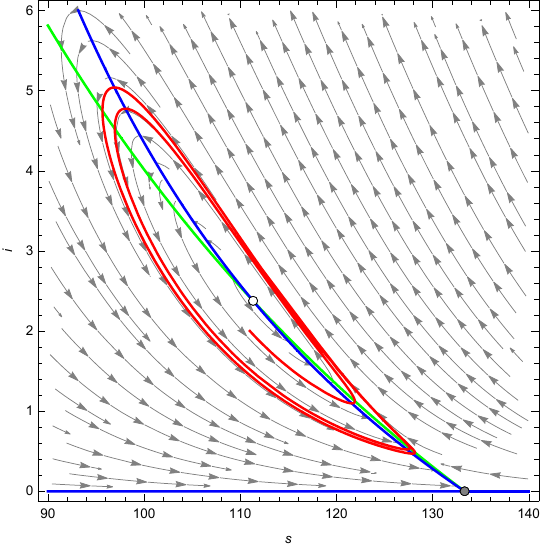}
        \caption{$(\s,\i)$- phase plot. The period of the cycle is $219.938 $ and the Floquet exponents are $\left\{-2.73386\times 10^{-8}, -0.0298991  \right\} $.}
    \end{subfigure}
    \caption{Time and phase-plot  corresponding  the neighboring region I  at the point $\left\{ \w=7.16058,\alpha=6.39679 \right\} $ in Figure \ref{f:fig6ns}.}
    \end{figure}

    {By continuity, cycles must still exist within region I, close to the boundary. This is somewhat surprising, since the DFE is
now stable}.

\sssec{Existence of Hopf-point at $B_1=\left\{ \w=5.15735 ,\alpha=4.60724 \right\}$}

This is  the   meeting point   of regions I,  II, III and VI. Here $\mR_0=1 \Lra $ $E_0=(133.333,0)=E_1$ is  a saddle point with eigenvalues $(-\mu=-0.12, 0)$, and the endemic point $E_2=(78.5251,  8.44639)$ is a {Hopf-point}  with eigenvalues  $\pm 0.152171\; Im  $, and $Tr[J(E_2)]=0$.
 An unstable limit cycle  around $E_2$  arises.
   
   \begin{figure}[H]
    \centering
    \begin{subfigure}[a]{0.4\textwidth}
    \centering
        \includegraphics[width=\textwidth]{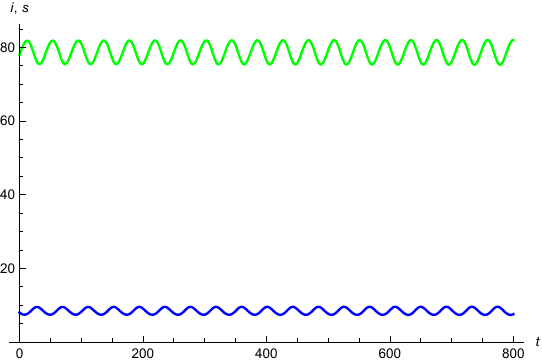}
        \caption{ $(s,i)$-time plot indicate the existence of a limit cycle surrounding $E_2$}
    \end{subfigure}%
    ~
     \begin{subfigure}[a]{0.38\textwidth}
     \centering
        \includegraphics[width=\textwidth]{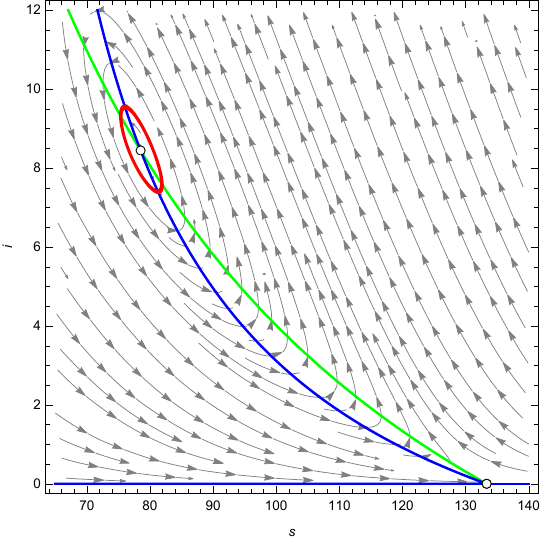}
        \caption{$(\s,\i)$- phase plot.  The cycle has  period  $42.2252$, and  Floquet exponents $\left\{0.00181877, -0.000205883 \right\}$}
    \end{subfigure}
  \caption{Time and phase-plot  corresponding to the boundaries $\mR_0=1, Tr[J(E_2)]=0$  at the point $B_1=\left\{ \w=5.15735 ,\alpha=4.60724 \right\}$ in Figure \ref{f:fig6ns}.}
    \end{figure}\la{f:limC}

\sssec{The point $B_2$ where regions II, III,VI,VIa meet}
At $B_2=(7.35966, 6.57463)$, one of the two solutions of $\mR_0=1, Tr(E_2)=0$, a stable limit cycle arises around the Hopf-point $E_2= (116.198,  1.77275) $ with eigenvalues $\pm 0.0391099\; Im$, while the other unstable fixed points $E_0=(133.333, 0)=E_1$ collide with eigenvalues $(-0.12, 0)$.

\begin{figure}[H]
    \centering
    \begin{subfigure}[a]{0.4\textwidth}
    \centering
        \includegraphics[width=\textwidth]{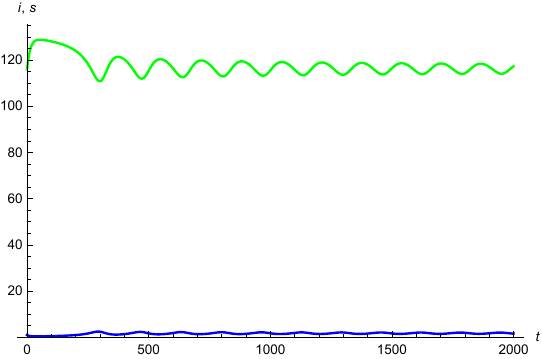}
        \caption{ $(s,i)$-time plot shows the existence of limit cycle surrounding $E_2$ }
    \end{subfigure}%
    ~
     \begin{subfigure}[a]{0.38\textwidth}
     \centering
        \includegraphics[width=\textwidth]{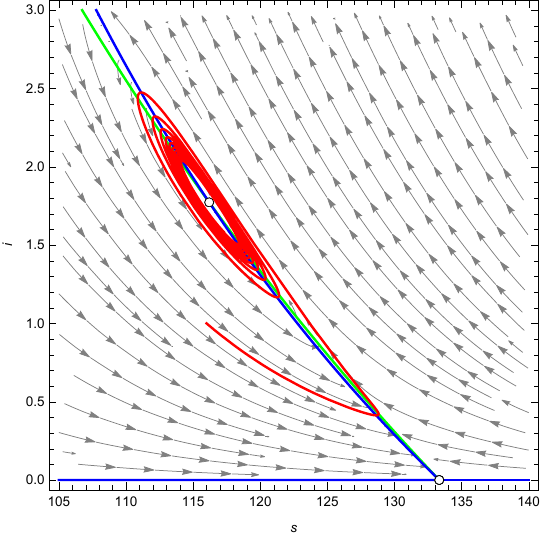}
        \caption{$(\s,\i)$- phase plot.  The cycle has period $160.805 $ and the Floquet exponents are $\left\{-0.0000365264, -0.000058307 \right\} $.}
    \end{subfigure}
    \caption{Time and phase-plot  corresponding   at the point $B_2=(7.35966, 6.57463)$ in Figure \ref{f:fig6ns}.}
    \end{figure}

\sssec{Phase-plot illustration  at the point $BT$ where regions III,V,VI,VIa meet}

Now moving to the point $BT=(\w=6.84183,\alpha= 6.20319)$ , the following figures indicate the convergence towards the \DFE\  $E_0=(133.333, 0)$ which is stable with eigenvalues $(-0.12, -0.0133232)$, and the remaining fixed points $E_1=E_2= (117.951, 1.56737) $ exhibits the existence of a Bogdanov-Takens bifurcation since the eigenvalues are $\pr{0, 0}$. 

\begin{figure}[H]
    \centering
    \begin{subfigure}[a]{0.4\textwidth}
    \centering
        \includegraphics[width=\textwidth]{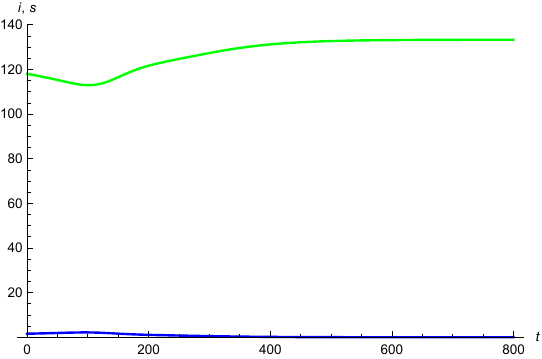}
        \caption{ $(s,i)$-time plot reveals convergence towards the \DFE }
    \end{subfigure}%
    ~
     \begin{subfigure}[a]{0.38\textwidth}
     \centering
        \includegraphics[width=\textwidth]{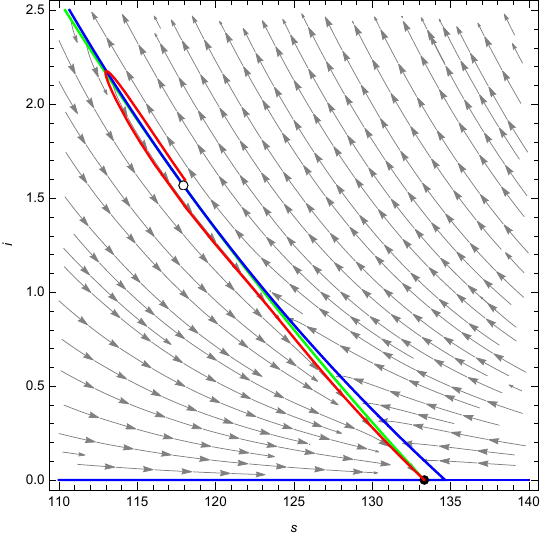}
        \caption{$(\s,\i)$- phase plot.  }
    \end{subfigure}
    \caption{Time and phase-plot   at the point $BT=(6.84183, 6.20319)$ in Figure \ref{f:fig6ns}.}
    \end{figure} 

\sssec{The point $H$ at the intersection of  $\mR_0=1$ and $\D=0$}

At $H=\left\{ \w=7.94466,\alpha=7.09723 \right\} $, the solution of $\mR_0=1$ and $\D=0$, the three fixed points coalesce as the figure below illustrates.

\begin{figure}[H]
    \centering
    \begin{subfigure}[a]{0.4\textwidth}
    \centering
        \includegraphics[width=\textwidth]{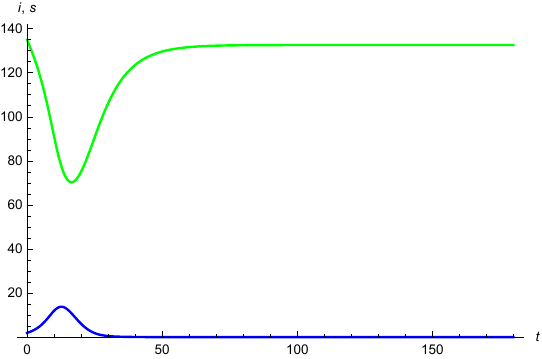}
        \caption{ $(s,i)$-time plot reveals convergence towards the \DFE }
    \end{subfigure}%
    ~
     \begin{subfigure}[a]{0.38\textwidth}
     \centering
        \includegraphics[width=\textwidth]{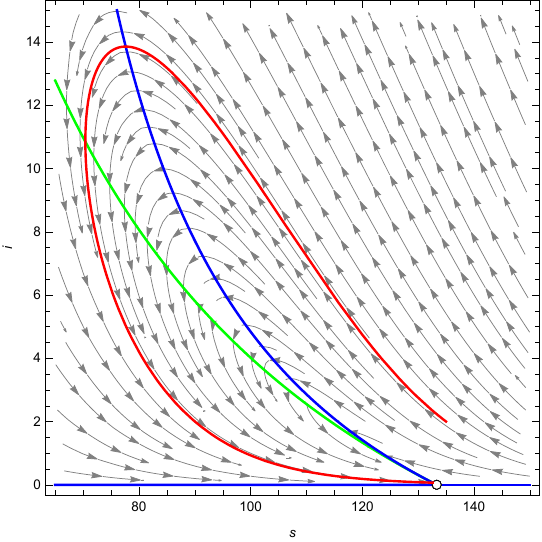}
        \caption{$(\s,\i)$- phase plot. The fixed points $E_0=E_1=E_2= (133.333, 0) $ collide with eigenvalues $ (-0.12, 0)$. }
    \end{subfigure}
    \caption{Time and phase-plot  corresponding to the boundary $\mR_0=1$ at the point $H$ in Figure \ref{f:fig6ns}, with $\left\{ \w=7.94466,\alpha=7.09723 \right\} \Lra  \mR_0=1 $.}
    \end{figure}

\sssec{Dynamical behaviour at region VIa}
Now moving inside the region VIa where we have bistability of $E_2=(122.3268072,1.080883)$ the \DFE\ $E_0=(133.333,0)$ with eigenvalues respectively,  $-0.01676672\pm 0.01328093\; Im  $ and $(-0.12, -0.0014)$.  The point $E_1=(130.4528,0.265)$ is a saddle point with eigenvalues $(-0.09477,0.0013)$.

\begin{figure}[H]
    \centering
    \begin{subfigure}[a]{0.4\textwidth}
    \centering
        \includegraphics[width=\textwidth]{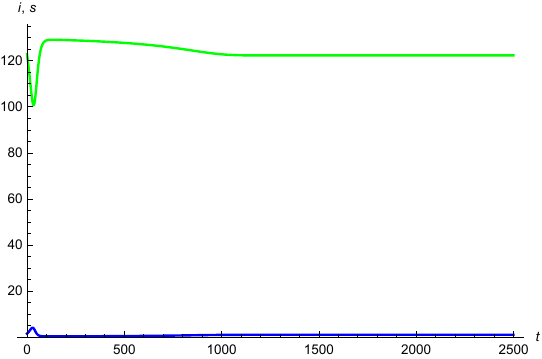}
        \caption{ $(s,i)$-time plot indicate the convergence towards $E_2$}
    \end{subfigure}%
    ~
     \begin{subfigure}[a]{0.38\textwidth}
     \centering
        \includegraphics[width=\textwidth]{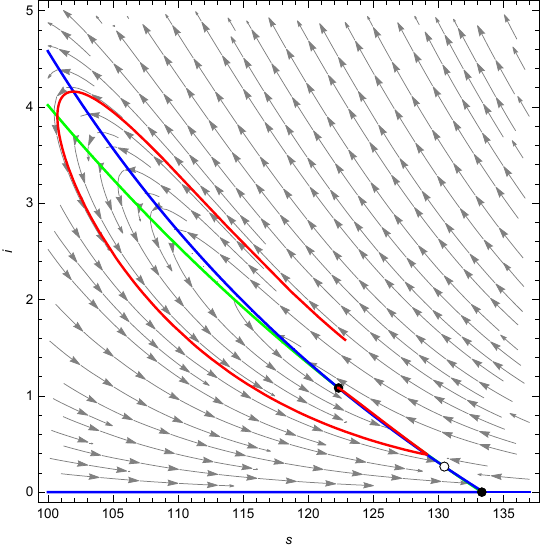}
        \caption{$(\s,\i)$- phase plot.  The stable attractor $E_2$ is surrounded by the red trajectory.}
    \end{subfigure}
    \caption{Time and phase-plot  corresponding to  region VIa  in Figure \ref{f:fig6ns},with $\left\{ \w=502/67,\alpha=41102/6131  \right\} $.}
    \end{figure}
    
 \sssec{The boundary $Tr(E_2)=0$ between  regions VIa  and I}
    At $T_2=\left\{ \w=2.93233,\alpha=3.69658  \right\}$, the stability of the DFE remains, and the endemic point The fixed point $E_2=(73.82, 9.770)$ becomes a Hopf point  with eigenvalues   $\pm 0.152 \; Im  $  and it's  surrounded by an unstable limit cycle. The DFE\  $E_0=(133.333,0)$ is stable with eigenvalues  $(-0.367, -0.12)$, and the point $E_1=(106.9, 2.971)$ is a saddle point with eigenvalues $(0.229, -0.0664)$.
    
    \begin{figure}[H]
    \centering
    \begin{subfigure}[a]{0.4\textwidth}
    \centering
        \includegraphics[width=\textwidth]{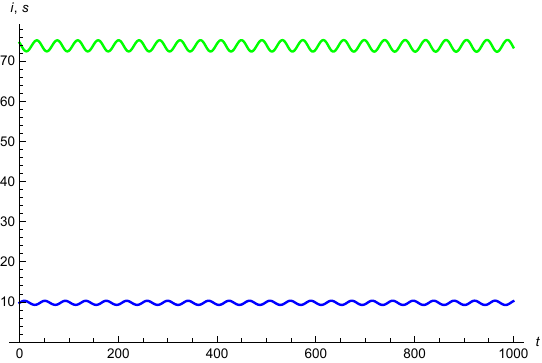}
        \caption{ $(s,i)$-time plot indicate the existence of a limit cycle surrounding $E_2$}
    \end{subfigure}%
    ~
     \begin{subfigure}[a]{0.38\textwidth}
     \centering
        \includegraphics[width=\textwidth]{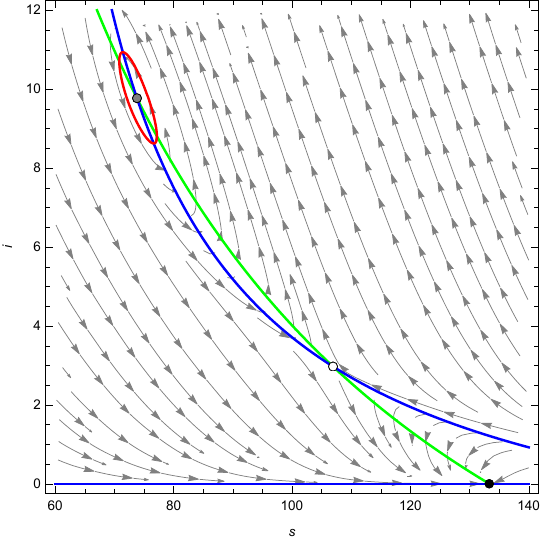}
        \caption{$(\s,\i)$- phase plot.  The cycle has period $41.7589$ and Floquet exponents $\left\{    0.00149016, 0.000030097 \right\}$. }
    \end{subfigure}
    \caption{Time and phase-plot  corresponding to  boundary $Tr(E_2)=0$ between  regions VIa  and I  in Figure \ref{f:fig6ns},with $\left\{ \w=2.93233,\alpha=3.69658  \right\} $.}
    \end{figure}

    \sssec{The boundary $\mR_0=1$ between  $B_2$  and $H$}
Moving on the boundary  $\mR_0=1$ towards the point $H=(7.94466,  7.09723)$ with $\w=7.91732,\; \alpha=7.0728$. The endemic point $E_2= (132.419,  0.0828753) $ is a stable with eigenvalues $(-0.111746, -0.0000342071)$. The \DFE\ $E_0=(133.333, 0)=E_1$ is a saddle point with eigenvalues $(-0.12, 0)$. The following figures show the convergence towards the \DFE.
\begin{figure}[H]
    \centering
    \begin{subfigure}[a]{0.4\textwidth}
    \centering
        \includegraphics[width=\textwidth]{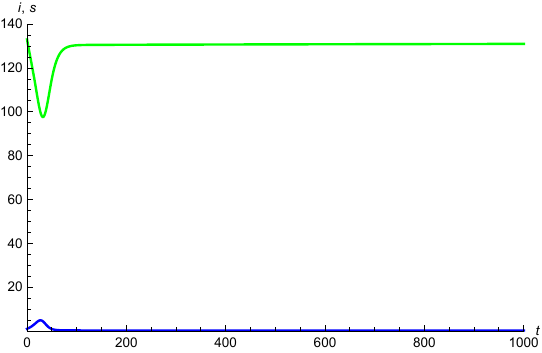}
        \caption{ $(s,i)$-time plot shows the convergence towards the \DFE }
    \end{subfigure}%
    ~
     \begin{subfigure}[a]{0.38\textwidth}
     \centering
        \includegraphics[width=\textwidth]{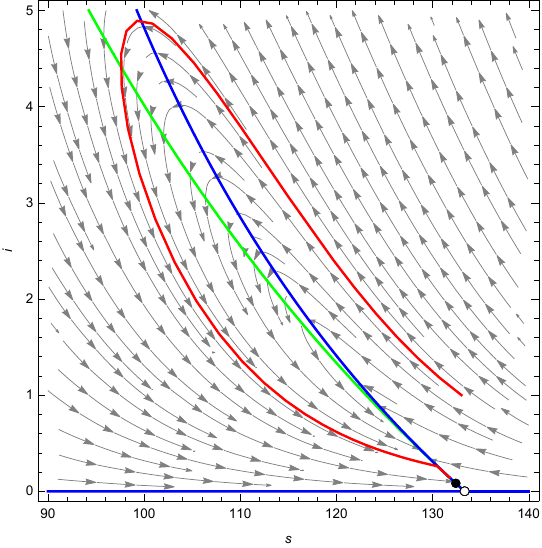}
        \caption{$(\s,\i)$- phase plot. }
    \end{subfigure}
     \caption{Time and phase-plot  corresponding to  boundary $\mR_0=1$ between $B_2$  and $H$  with $\w=7.91732,\; \alpha=7.0728$ in Figure \ref{f:fig6ns}.}
    \end{figure}

   \sssec{The boundary $\mR_0=1$ between regions III and VI}\la{s:IItoI}
 At the point $R_1 =\pr{\w=4.62051, \alpha=4.12766}$, and moving downwards from $B_1=(5.15735,  4.60724)$, the following figures show the convergence towards the endemic point $E_2=(72.8175,  10.0732)$ which is a stable attractor  with eigenvalues  $-0.017169 \pm 0.173632 \; Im  $. The fixed points  $E_0=E_1=(133.333,0)$ collide  into a saddle with eigenvalues   $(-0.12,0)$. 
  
   \begin{figure}[H]
    \centering
    \begin{subfigure}[a]{0.4\textwidth}
    \centering
        \includegraphics[width=\textwidth]{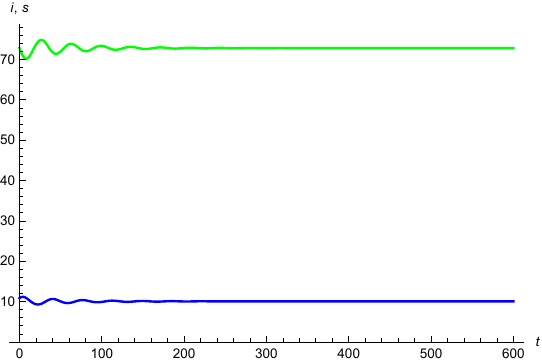}
        \caption{ $(s,i)$-time plot indicate the convergence towards $E_2$}
    \end{subfigure}%
    ~
     \begin{subfigure}[a]{0.38\textwidth}
     \centering
        \includegraphics[width=\textwidth]{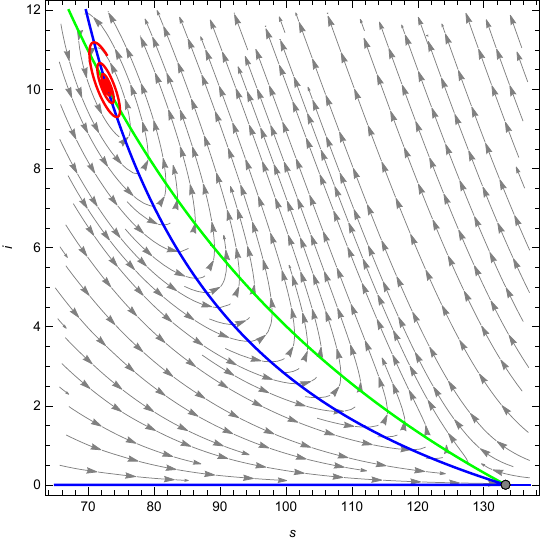}
        \caption{$(\s,\i)$- phase plot. }
    \end{subfigure}
    \caption{Time and phase-plot  corresponding to  boundary $\mR_0=1$ between regions III and VI  at the point $R_1 =\pr{\w=4.62051, \alpha=4.12766}$ in Figure \ref{f:fig6ns}.\label{f:R1}}
    \end{figure}

Starting to move upwards on the boundary closer to the ``B point" $B_1$ where regions $I,II,III,VI$ meet  reveals in next Figure \ref{f:hom} a stable limit cycle, surrounding the unique unstable interior point $E_2=(79.9285,  8.0827)$ with eigenvalues  $0.00347509 \pm 0.146926 \; Im  $, which is ``almost  homoclinic \wrt\   $E_0=E_1=(133.333,0)$"   (with eigenvalues  $(-0.12,0)$).

    \begin{figure}[H]
    \centering
    \begin{subfigure}[a]{0.4\textwidth}
    \centering
        \includegraphics[width=\textwidth]{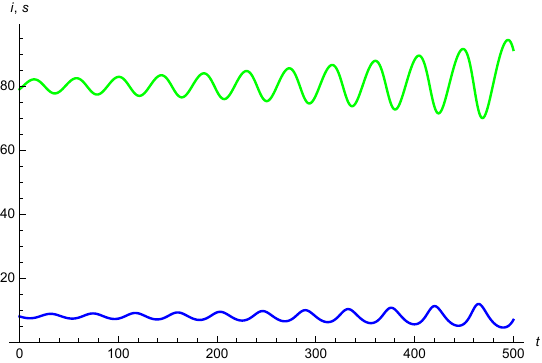}
        \caption{ $(s,i)$-time plot indicate the existence of a limit cycle surrounding $E_2$}
    \end{subfigure}%
    ~
     \begin{subfigure}[a]{0.38\textwidth}
     \centering
        \includegraphics[width=\textwidth]{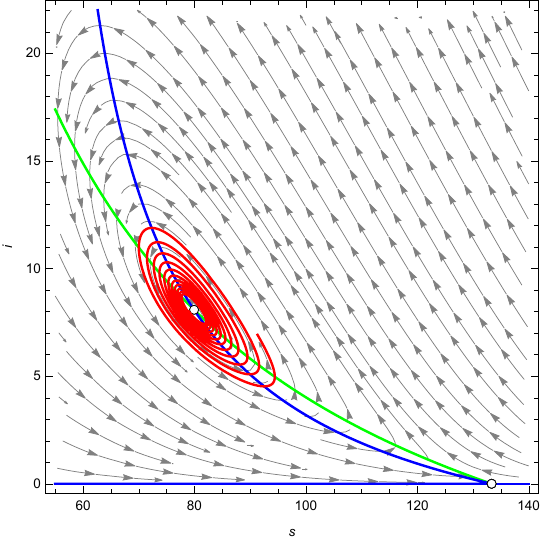}
        \caption{$(\s,\i)$- phase plot.  The   period is $40.1628$, and the Floquet exponents $-0.00353727 \pm 0.000583175\; Im$}
    \end{subfigure}
    \caption{Time and phase-plot  corresponding to  boundary $\mR_0=1$ between regions I and II  in Figure \ref{f:fig6ns}, with $\left\{ \w=  723/137  ,\alpha= 16147/3425 \right\} \Lra  \mR_0=1$.}
    \end{figure}

 We  continue now moving further upwards from  $B_1=\left\{ \w=4.60724,\alpha=5.15735  \right\}$ and getting closer to $BT=(6.84183,  6.20319)$, on the boundary $\mR_0=1$ between regions II and I with $\left\{ \w=25/4,\alpha=67/12  \right\}$.
Since $\mR_0=1$, the DFE equals $E_1=(133.333,0)$ with eigenvalues  $(0,-0.12)$ , and the following figure reveals {an unstable, ``almost homoclinic" limit cycle}, surrounding  the endemic point $E_2=(93.5174, 5.13535)$ with eigenvalues  $0.0226741 \pm 0.0987278 \; Im  $.
 
\begin{figure}[H]
    \centering
    \begin{subfigure}[a]{0.4\textwidth}
    \centering
        \includegraphics[width=\textwidth]{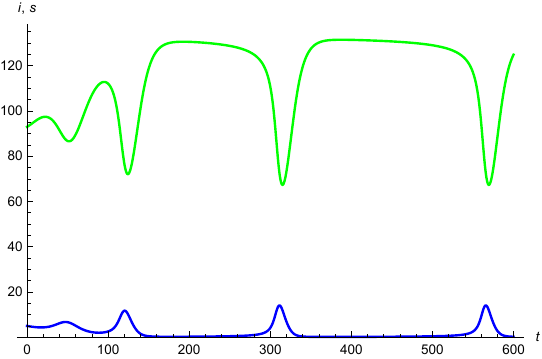}
        \caption{ $(s,i)$-time plot indicate the existence of a limit cycle surrounding $E_2$}
    \end{subfigure}%
    ~
     \begin{subfigure}[a]{0.38\textwidth}
     \centering
        \includegraphics[width=\textwidth]{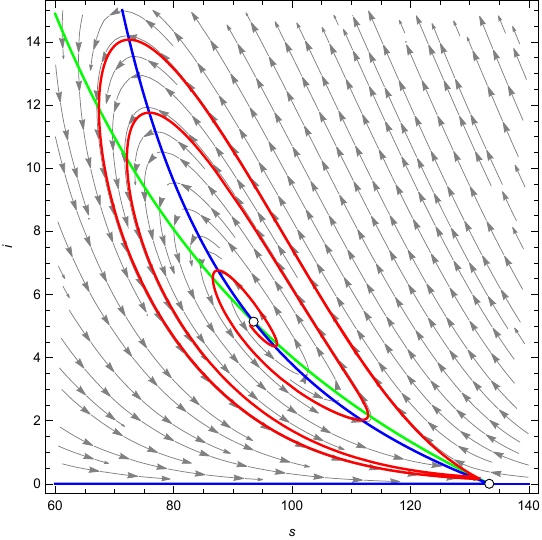}
        \caption{$(\s,\i)$- phase plot. The   period is $254.473$, and the Floquet exponents $\left\{ -4.59451\times10^{-8}, -0.0659782\right\}$}
    \end{subfigure}
    \caption{Time and phase-plot  corresponding to  boundary $\mR_0=1$ between regions I and II  in Figure \ref{f:fig6ns}.\label{f:hom}}
    \end{figure}

\sssec{The bistability region VI (with DFE and $E_2$ stable and $E_1$  saddle)}

In this region, at the point $Q_{VI}=\pr{\w=0.078125,\alpha=0.15625}$, both $E_2=(45.712,23.542)$ and the DFE $E_0=(133.33,0)$ are locally stable with eigenvalues $\pr{-0.176845\pm 0.265287\; Im}$ and $(-1.1,-0.12)$, respectively. The fixed point $E_1=(132.24,0.098)$ is a saddle point with the corresponding eigenvalues  $(0.4908,-0.1188)$.
 We illustrate in the following figures the convegence towards this two points and the corresponding phase-plot.

\begin{figure}[H]
    \centering
    \begin{subfigure}[a]{0.4\textwidth}
        \includegraphics[width=\textwidth]{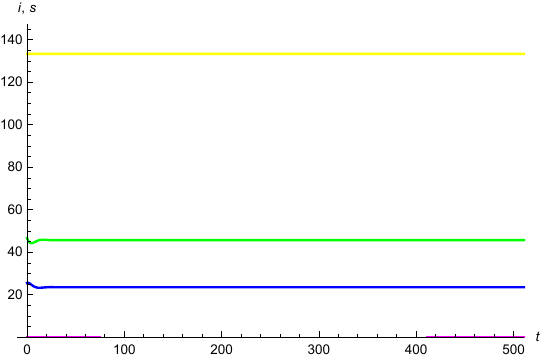}
        \caption{$(\s,\i)-$time plots of indicate convergence  towards $E_2$ and DFE.}
        \label{fig:sol6}
    \end{subfigure}
    ~
    \begin{subfigure}[a]{0.38\textwidth}
        \includegraphics[width=\textwidth]{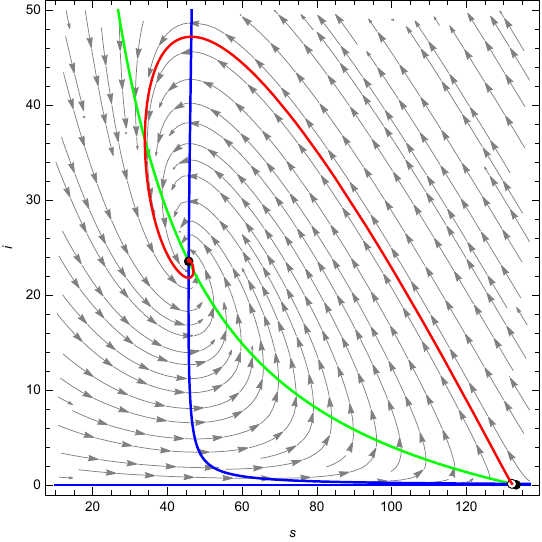}
        \caption{$(\s,\i)$- phase plot.}
        \label{fig:sp6}
    \end{subfigure}
    \caption{Time and phase-plot corresponding to region VI at $Q_{VI}=\pr{\w=0.078125,\alpha=0.15625}$ in figure \ref{f:fig6ns}. }
    \end{figure}

{We have not found
the claimed  homoclinic cycle in  \cite[Fig. 1c)]{Gupta}}


{The problem of whether oscillations are possible in region VI away from the $B_1$ point is still open.}

{We skip the region I, where the DFE is the only stable point}.

\section{ Appendix: A description of  the accompanying Mathematica notebook}\la{s:num}
Writing a notebook for studying symbolically and numerically a dynamical system with seven parameters is a highly non-trivial task. This may be the explanation why a very small percentage of the literature is accompanied by notebooks.

  We choose to organise our notebook \cite{NotM1} as follows:
\begin{itemize}
 \item First cell enumerates the parameters to be studied.
 \item If we may want to use other parameters at a later time,
  in a future extension, we provide a condition setting them equal to $0$. We include here also notable particular cases.
\item  Afterwards, we provide definitions of  simple fundamental quantities, which may be easily  computed  without Mathematica, and are ``correct beyond any doubt".
  In our case, we define the \brn\ $\mR_0$, and some  critical parameters which make it equal to $1$ ($\b,\eta$). In this way, we can specify the fundamental hypersurface $\mR_0=1$ by fixing either $\b$ or $\eta$.
\item   We provide then the essential conditions on the parameters (positivity).

\item  There will be also fundamental quantities which require a  long time for computing and simplifying symbolically; these have been saved in the ``package" "def.m"  and are loaded next.

In our cases, they are the trace at the second endemic point, which is long that it requires 2 pages to display, and the discriminant $\Delta$. This equation reveals that the discriminant is only defined for $\w v_1 < \mu$, and that the curve $\Delta=0$ may be explicitized with respect to $\a$, which  is clearly useful.

\item We continue with  some numerical conditions used in first tests, and with some conditions for switching between the parameters, for displaying alternative representations of the answers.
\item  The second cell contains the definition of the model, and fundamental quantities (fixed points, Jacobians, etc).
Here we can either produce our own code, or use EcoEvo. We will provide both versions, since  each has its advantages.

 \end{itemize}

\section*{Declarations}

\textbf{Conflict of Interest:} The authors have no competing interests to declare that are relevant to the content of this article.
 
 \bibliographystyle{amsalpha}

\bibliography{Pare38}

\end{document}